\documentclass[a4paper,11pt]{article}
\usepackage{amssymb,amsmath,amsfonts,epsf,amsmath,tikz}
\usepackage{enumerate}
\usetikzlibrary{decorations.pathreplacing,automata,calc,positioning}
\usepackage{framed}

\newtheorem{theorem}{\bf Theorem}[section]

\newtheorem{lemma}[theorem]{\bf Lemma}

\newtheorem{proposition}[theorem]{\bf Proposition}
\newtheorem{conjecture}[theorem]{\bf Conjecture}

\newtheorem{definition}[theorem]{\bf Definition}
\newtheorem{observation}[theorem]{\bf Observation}

\newcommand{\D}{{Dominator}}
\newcommand{\St}{{Staller}}

\renewcommand{\ggg}{{\gamma_g}}
\newcommand{\gggp}{{\gamma_g'}}

\newcommand{\efface}[1]{}
\newcommand{\proof}{\noindent{\bf Proof.\ }}
\newcommand{\qed}{\hfill $\square$ \medskip}

\newcommand{\go}{D-game}
\newcommand{\gt}{S-game}
\newcommand{\sff}{weakly $S(K_{1,3})$-free}

\begin{document}

\title{The 3/5-conjecture for \sff\ forests}

\author{Simon Schmidt $^{a}$}

\date{}

\maketitle

\begin{center}
$^a$ Institut Fourier, SFR Maths \`a Modeler, Joseph Fourier's University

\noindent
100 rue des Maths, BP 74

\noindent
38402 St Martin d'H\`eres, France

\noindent
simon.schmidt@ujf-grenoble.fr
\end{center}


\begin{abstract}
The $3/5$-conjecture for the domination game states that the game domination numbers of an isolate-free graph $G$ on $n$ vertices are bounded as follows: $\ggg(G)\leq  \frac{3n}5 $ and $\gggp(G)\leq  \frac{3n+2}5 $. 
Recent progress have been done on the subject and the conjecture is now proved for graphs with minimum degree at least $2$. 
One powerful tool, introduced by Bujt\'as is the so-called greedy strategy for \D. In particular, using this strategy, she has proved the conjecture for isolate-free forests without leafs at distance $4$. In this paper, we improve this strategy to extend the result to the larger class of \sff\ forests, where a \sff\ forest $F$ is an isolate-free forest without induced $S(K_{1,3})$, whose leafs are leafs of $F$ as well.
\end{abstract}

\noindent{\bf Keywords:} domination game; $3/5$-conjecture

\medskip
\noindent{\bf AMS Subj. Class.:} 05C57, 91A43, 05C69

\section{Introduction}
The domination game, introduced five years ago in \cite{bresar-2010} is played on an arbitrary graph $G$ by two players,
{\em Dominator} and {\em Staller}. They alternately choose a vertex from $G$ such that at least one previously undominated vertex becomes dominated. 
The game ends when no move is possible. \D\ aims to end the game as soon as possible, while \St\ wants to prolong it. By {\em \go\ } (resp. {\em \gt}) we mean a game in which Dominator (resp. Staller)
plays first. Assuming that both players play optimally, the {\em \go\ domination number} $\gamma_g(G)$ (resp. the {\em \gt\ domination number} $\gamma_g'(G)$) of a graph $G$, denotes the total number of chosen vertices during \go\ (resp. \gt) on $G$.
There is already a flourishing amount of works on the subject, see for various examples \cite{kosmrlj-2014,dorbec-2015,complexity-2014+,brdo-2014,dorbec-2015+,javad-2015+}. The most outstanding conjecture concerning this game is probably the so-called $3/5$-conjecture. 
(For related developments concerning the total domination game see \cite{henning-2015,henning-2015+}.)
\begin{conjecture}\textbf{(Kinnersley, West and Zamani, \cite{kinnersley-2013})} If $G$ is an isolate-free graph on $n$ vertices, then $$\ggg(G)\leq \frac{3n}5 \text{ and } \gggp(G)\leq\frac{3n+2}5.$$ 
\end{conjecture}
If true, the bounds are known to be tight, even if we restrict ourself to forests. See \cite{bresar-2013} for a study of the structure of graphs reaching these bounds. 
Recently Bujt\'as has made a breakthrough by introducing a powerful greedy strategy for \D\ \cite{bujtas-2013}. In particular, using this technique, she has proved the conjecture for a large subclass of forests \cite{bujtas-2013,bujtas-2014}. In that case, the upper bound for \gt\ turns two be even better.
\begin{theorem}[Bujt\'as 2014]
 If $F$ is an isolate-free forest on $n$ vertices, without leafs at distance $4$, then
$$\ggg(F)\leq \frac{3n}5 \text{ and } \gggp(F)\leq\frac{3n+1}5.$$
\end{theorem}
This greedy strategy could also be applied to graphs with minimum degree at least $3$. 
In \cite{bujtas-2015}, the conjecture is proved for this class of graph. Even more recently, Henning and Kinnersley~\cite{henning-2014+} established the truth of the $3/5$-conjecture over the class of graphs with minimum degree at least~$2$. Hence the $3/5$-conjecture remains open only for graphs with pendant vertices. 
In this paper, we prove the $3/5$ conjecture for a larger class of forests, which contains the forests without leafs at distance 4. Even if forests deserve interest by themselves, we emphasize that 
the general conjecture cannot be  easily reduced to the one on forests. There are actually graphs which have greater domination number than any of their spanning trees \cite{brkl-2013b}.  

\medskip
Let $S(K_{1,3})$ be the graph obtained by subdividing once all the edges of the star $K_{1,3}$ (see Figure~\ref{fig:SK}). A {\em \sff\ forest} is a forest which does not contain:
\begin{itemize}
 \item isolated vertices and
 \item induced $S(K_{1,3})$ whose leafs are leafs in the original forest. 
\end{itemize}
\begin{figure}[!ht]
\begin{center}
\begin{tikzpicture}[scale=1,style=thick]
\def\vr{3pt}
\draw(0,0) coordinate(x1)--(0.5,0.5) coordinate (x2)--(1,1) coordinate(x3)--(1.5,0.5) coordinate (x4)--(2,0) coordinate (x5);
\draw(1,1.7) coordinate(x6)--(1,2.4) coordinate (x7)--(x3);

\draw(x1)[fill=white] circle(\vr);\draw(x2)[fill=white] circle(\vr);
\draw(x3)[fill=white] circle(\vr);\draw(x4)[fill=white] circle(\vr);
\draw(x5)[fill=white] circle(\vr);\draw(x6)[fill=white] circle(\vr);
\draw(x7)[fill=white] circle(\vr);
\end{tikzpicture}
\end{center}
\caption{$S(K_{1,3})$}
\label{fig:SK}
\end{figure}
The class of \sff\ forests obviously contains the isolate-free forests without leafs at distances $4$. For this class of forests we are able to extend Bujt\'as greedy strategy to prove the following more general theorem.

\begin{theorem}
If $F$ is a \sff\ forest on $n$ vertices, then $$\displaystyle \ggg(F)\leq\frac{3n}{5} \text{ and } \gggp(F)\leq\frac{3n+2}{5}.$$
\end{theorem}
The upper bounds are tight in both cases. For example, by \cite[Theorem 3.7]{bresar-2013}, trees from Figure~\ref{fig:ext} (a) reach the bound for \go. For \gt, we easily verify that forests from Figure~\ref{fig:ext} (b) attain the bound as well.  
\begin{figure}[!ht]
\begin{center}
\begin{tikzpicture}[scale=1,style=thick]
\def\vr{3pt}
\draw(0,0) coordinate(x1)--(0.5,0.25) coordinate (x2)--(1,0.5) coordinate(x3)--(1.5,0.25) coordinate (x4)--(2,0) coordinate (x5);
\draw(0,2) coordinate(x6)--(0.5,2.25) coordinate (x7)--(1,2.5) coordinate(x8)--(1.5,2.25) coordinate (x9)--(2,2) coordinate (x10);
\draw(0,2.8) coordinate(x11)--(0.5,3.05) coordinate (x12)--(1,3.3) coordinate(x13)--(1.5,3.05) coordinate (x14)--(2,2.8) coordinate (x15);
\draw(x8)--(x13);\draw[dashed](x8)--(x3);

\draw(x1)[fill=white] circle(\vr);\draw(x2)[fill=white] circle(\vr);
\draw(x3)[fill=white] circle(\vr);\draw(x4)[fill=white] circle(\vr);
\draw(x5)[fill=white] circle(\vr);\draw(x6)[fill=white] circle(\vr);
\draw(x7)[fill=white] circle(\vr);\draw(x8)[fill=white] circle(\vr);
\draw(x9)[fill=white] circle(\vr);\draw(x10)[fill=white] circle(\vr);
\draw(x11)[fill=white] circle(\vr);\draw(x12)[fill=white] circle(\vr);
\draw(x13)[fill=white] circle(\vr);\draw(x14)[fill=white] circle(\vr);
\draw(x15)[fill=white] circle(\vr);

\draw(x3)[above right] node {$x_n$};
\draw(x13)[above right] node {$x_1$};
\draw(x8)[above right] node {$x_2$};
\draw(1,0)[below] node {(a)};

\draw(4,0) coordinate(y1)--(4.5,0.25) coordinate (y2)--(5,0.5) coordinate(y3)--(5.5,0.25) coordinate (y4)--(6,0) coordinate (y5);
\draw(4,2) coordinate(y6)--(4.5,2.25) coordinate (y7)--(5,2.5) coordinate(y8)--(5.5,2.25) coordinate (y9)--(6,2) coordinate (y10);
\draw(4,2.8) coordinate(y11)--(4.5,3.05) coordinate (y12)--(5,3.3) coordinate(y13)--(5.5,3.05) coordinate (y14)--(6,2.8) coordinate (y15);
\draw(y13)--(5,4.1) coordinate (y0);

\draw(y1)[fill=white] circle(\vr);\draw(y2)[fill=white] circle(\vr);
\draw(y3)[fill=white] circle(\vr);\draw(y4)[fill=white] circle(\vr);
\draw(y5)[fill=white] circle(\vr);\draw(y6)[fill=white] circle(\vr);
\draw(y7)[fill=white] circle(\vr);\draw(y8)[fill=white] circle(\vr);
\draw(y9)[fill=white] circle(\vr);\draw(y10)[fill=white] circle(\vr);
\draw(y11)[fill=white] circle(\vr);\draw(y12)[fill=white] circle(\vr);
\draw(y13)[fill=white] circle(\vr);\draw(y14)[fill=white] circle(\vr);
\draw(y15)[fill=white] circle(\vr);\draw(y0)[fill=white] circle(\vr);

\draw(y3)[above right] node {$y_n$};
\draw(y13)[above right] node {$y_1$};
\draw(y8)[above right] node {$y_2$};
\draw(y0)[above right] node {$y_0$};
\draw(5,0)[below] node {(b)};

\end{tikzpicture}
\end{center}
\caption{Extremal \sff\ forests}
\label{fig:ext}
\end{figure}

In the second section of the article, we are going to prove the upper bounds of our main theorem. 
The first subsection describes Bujt\'as' greedy strategy and introduces an improvement of it, the so-called breaking-$P_5$ strategy. The second subsection deals with the proof 
of the main theorem for \go. Finally, in the last subsection, we explain briefly how to modify the proof to extend the result to \gt.

\section{Proof of the main theorem}

We start by proving the upper bound for \go. All the results of the two following subsections are stated assuming that \go\ is played. 
In the last subsection, we explain how to modify the proof to get the desired upper bound for \gt.

\subsection{The greedy strategy and the breaking-$P_5$ strategy}\label{sec:defstrat}

First we recall Bujt\'as' greedy strategy for \D\ and give definitions and background needed to understand it. 
Our strategy will be compatible with the greedy one. In fact, it is only a refinement of the first phase of this strategy.
Let $F$ be an isolate-free forest on $n$ vertices. All vertices will be colored with one of the following three colors: white ($W$), blue ($B$) and red ($R$). In addition to colors, each vertex of $F$ has also a value. 
It is this value which ables to define a greedy like strategy. The coloring and the values will change all along the game. At the end of turn $k\geq 0$, the color and value of a vertex are defined as follows.
\begin{itemize}
 \item A vertex is white and its value is $3$ if it is not already dominated.
 \item A vertex is blue and its value is $2$ if it is dominated, but at least one of its neighbors is not.
 \item A vertex is red and its value is $0$ if $N[v]$ is entirely dominated.
 \end{itemize}

Note that, when a vertex is selected by a player, it is always going to turn red during the current turn. We denote respectively by $W_k(F)$, $B_k(F)$ and $R_k(F)$, the set of white, blue and red vertices at the end of turn $k\geq 0$. 
When $F$ will be clear from the context, we will omit it in the notations. At the end of turn $k\geq 0$, we define the residual forest $F_k$ as the forest obtained from $F$ by deleting all the red vertices and all the edges between two blue vertices (note that $F_0=F$).
Since all the legal moves are either white or blue vertices and all the blue vertices are already dominated, the game after $k$ moves is the same played on the forest $F$ or played on the residual forest $F_k$.
The following statement proved in \cite{bujtas-2014} is particularly useful. In particular, the first item ensures that any white leaf in a residual forest is also a leaf in the original forest. In other words, no white leaf is created during the game. 
\begin{lemma}\label{lem:W-leaf}
Let $F$ be an isolate-free forest and $k\geq 0$.
 \begin{enumerate}
  \item If $v$ belongs to $W_k$, then $v$ has the same neighborhood in $F_k$ as in $F$.
  \item If $F$ has no isolated vertex, then neither has $F_k$. 
 \end{enumerate}
\end{lemma}

For $k\geq 0$, the value of the residual forest $F_k$ is the sum of the values of all its vertices. We denote this value by $p(F_k)$. By definition, $p(F)=3n$. At each turn, the value of the residual forest decreases. We say that a player seizes $s$ points during turn $k\geq1$, when $p(F_{k-1})-p(F_k)=s$. 
In the greedy strategy the game will be divided into four phases described below. \D\ will always start by applying, if possible, Phase 1 of the strategy. 
For $i\in\{2,3,4\}$, Phase $i$ starts at \D's turn only if Phase $i-1$ is no longer possible. If at this moment, the strategy of Phase $i$ is not applicable, then it will be skipped. When \D\ 
starts to play in a new phase, he will never go back to a previous one, even if a change in the residual graph might cause this phase applicable again.

\medskip\noindent
\textbf{Greedy strategy.}

\begin{itemize}
 \item \textbf{Phase 1.} At his turn \D\ seizes at least $7$ points and at least two vertices turn red during this turn.
 \item \textbf{Phase 2.} At his turn \D\ seizes at least $7$ points.
 \item \textbf{Phase 3.} At his turn \D\ seizes at least $6$ points. Moreover, he applies the following two rules.
	  \begin{itemize}
	      \item  \D\ selects a vertex which ensures the maximal possible gain achievable at this turn.
	      \item  Under the above rule, \D\ always prefers to play a white stem in the residual forest, which has a white leaf neighbor.
	   \end{itemize}
 \item \textbf{Phase 4.} At his turn \D\ seizes at least $3$ points.
\end{itemize}
Roughly speaking, Bujt\'as' proof is to show that during a turn of Dominator and the turn of \St\ which follows, the value of the residual forest is in average decreased by $10$ points.
It implies directly that after at most $3n/5$ moves, the value of the residual graph is $0$, which means the game is over. But it works only if we do not deal with what Bujt\'as called 
critical $P_5$. 

\begin{definition}
 A critical $P_5$ in $F_k$, with $k\geq 0$ is a path of length $4$ such that the following holds. 
 \begin{itemize}
  \item All the vertices are in $W_k$, save the center of the path which is in $B_k$.
  \item Both ends of the path are leafs in $F_k$. 
 \end{itemize}
\end{definition}
The blue center of a critical $P_5$ is called a {\em critical center}. 
It is proved in \cite{bujtas-2014}, that these critical centers could only appear during Phase 1. Since no new white leaf is created along the game, the following observations are also easy to prove.
\begin{observation}\label{ob:1}
 If a vertex in $B_k$ is not a critical center at the end of turn $k\geq1$, then it is not a critical center at the end of any turn $k'\geq k$.
\end{observation}
In order to state precisely the results obtained by Bujt\'as, we need more definitions related to this phase.
Let $k^*\geq 0$ be the turn at the end of which Phase 1 is over. If Phase 1 is skipped, it is easy to show that the $3/5$-conjecture holds (the forest is actually a disjoint union of edges). For, we assume in the rest of the paper that $k^*\geq 1$.
The number of critical centers at the end of turn $k^*$ will be denoted by $c^*$.
For each turn of Phase 1, when Dominator decreases the value of the graph by more than $7$, he does better than what it is expected during this phase. It is the same for \St\ when she 
decreases the value by more than $3$ points. These bonus points are defined as follows.
\begin{definition}
 For $k\in\{1,...,k^*\}$, set:
$$e_k=\begin {cases}
       p(F_{k-1})-p(F_k)-7&\text{if turn $k$ is due to \D,}\\
       p(F_{k-1})-p(F_k)-3&\text{if it is due to \St.}
      \end {cases}$$
\end{definition}
The sum of all the bonus points earned during Phase 1 will be denoted by $e^*$, that is $e^*=e_1+\cdots+e_{k^*}$. In \cite[Theorem~2]{bujtas-2014}, Bujt\'as proved the following result.

\begin{proposition} \label{prop:buj}
 Let $F$ be an isolate-free forest on $n$ vertices. If \D\ plays \go\ according to the greedy strategy, then the number of turns will be at most $\displaystyle \displaystyle \frac{3n-e^{*}+c^{*}}5 $.
\end{proposition}
By Lemma~\ref{lem:W-leaf}, no new white leaf is created all along the game. It follows that for an isolate-free forest $F$ without leafs at distance $4$, $c^*$ is always equal to $0$. 
Hence, the above proposition shows that the $3/5$-conjecture holds for \go\ when restricted to this class of forests. For \gt\, similar arguments show that for this restricted class the upper 
bound is slightly better than the one announced by the conjecture: $\gggp(F)\leq \frac{3|V(F)|+1}2$.

\medskip
We now improve the greedy strategy in order to apply it to the larger class of \sff\ forests.
Let $F$ be an isolate-free forest on $n$ vertices. From now on, we assume that $F$ is arbitrarily rooted. 
More precisely, if $F$ is a forest with $m\geq 1$ connected components $T_1$,...,$T_m$, then for any $i\in\{1,...,m\}$ we select some $r_i\in V(T_i)$ to be the root of the tree $T_i$. 
For every vertex $v\in V(T_i)$, the {\em height of $v$} is defined by $h(v)=d(r_i,v)$.

\medskip\noindent
We define $\mathcal P (F)$ as the set of paths of length $4$ between two leafs of $F$. We say that a vertex $d\in V(F)$ is a {\em dangerous center after $k\geq 0$ moves}, if the followings holds.
\begin{itemize}
 \item The vertex $d$ is the center of a path $P$ in $\mathcal P(F)$.
 \item All the vertices of $P$ are in $W_k$.
\end{itemize}
The set of dangerous centers after $k\geq 0$ moves will be denoted by $D_k(F)$. When $F$ will be clear from the context, we will only write $D_k$. 
Because no new white leaf is created in the residual forests along the game (see Lemma ~\ref{lem:W-leaf}), the two following results clearly hold.%
\begin{lemma}
If $k'\geq k\geq 0$ then $D_{k'} \subseteq D_k$, that is no new dangerous center is created during the game.
\end{lemma}
\begin{lemma}\label{lem:0}
 If a vertex in $W_k$ does not belong to $D_k$, with $k\geq 0$, then it could not be a critical center at the end of any turn $k'> k$.
\end{lemma} 
The following lemma is important, because it gives an alternative definition of the class of \sff\ forests. 
\begin{lemma}\label{lem:a1a2}
If $F$ is a \sff\ forest and $d$ is a dangerous center in $D_0(F)$, then there are exactly two vertices $a_1,a_2\in N(d)$, such that any path in $\mathcal P(F)$ whose center is $d$ contains $a_1$ and $a_2$. 
\end{lemma}
\proof Let $l_1a_1da_2l_2$ be a path in $\mathcal P(F)$ whose center is the dangerous center $d$. Let $f$ be a leaf of a path $P$ in $\mathcal P(F)$ whose center is also $d$. 
We have to prove that the unique neighbor $u$ of $f$ is either $a_1$ or $a_2$. By way of contradiction, assume that $P$ contains  neither $a_1$ nor $a_2$. 
Then the subtree induced by $\{l_1,l_2,f,a_1,a_2,u\}$ is clearly isomorphic to $S(K_{1,3})$. Since $l_1$, $l_2$ and $f$ are leafs of $F$, the forest is not \sff. A contradiction.
\qed

At least one of the two vertices $a_1$, $a_2$ defined in Lemma~\ref{lem:a1a2} has a height strictly greater than $d$. 
Without lost of generality, we always assume that $h(a_2)>h(d)$ and we say that $a_2$ is the {\em $P_5$-child} of $d$. We introduce now the breaking-$P_5$ strategy. The notation $a_1$ and $a_2$ will always refer to 
the two vertices defined in the above lemma. We will say that $a_1$ is the vertex {\em related} to the dangerous center $d$.  

\medskip\noindent
\textbf{Breaking-$P_5$ strategy.}
\begin{itemize}
 \item \textbf{Phase 1.0:} if $D_{k-1}$ is not empty at his turn $k\geq 1$, \D\ selects a vertex $d\in D_{k-1}$ with maximum height and he plays the vertex $a_1$ related to $d$. To emphasize that $d$ is not the vertex played
 during this turn, we will always say that $d$ is {\em elected}.
 \item \textbf{Phase 1.1:} at his turn, \D\ gets at least 7 points and at least two vertices turn red.
 \item \textbf{Phases 2 to 4:} these phases are the same as for the greedy strategy.
 \end{itemize}
The phases follow one another in the same way as for the greedy strategy.
We write $k^{**}\geq 0$ for the number of turns in Phase 1.1, $c^{**}$ for the number of critical centers at the end of turn $k^{**}$ and $e^{**}=e_1+\cdots+e_{k^{**}}$. 
If $k^{**}=0$, then we are in the settings of Bujt\'as' proof. Hence, we assume now that $k^{**}\geq 1$.
As we will see in the following proposition, the breaking-$P_5$ strategy is also a greedy strategy and we have a result similar to Proposition~\ref{prop:buj}.

\begin{proposition}\label{prop:1}
 Let $F$ be an isolate-free forest on $n$ vertices. If \D\ plays \go\ according to the breaking-$P_5$ strategy, then the number of turns will be at most $\displaystyle \frac{3n-e^{**}+c^{**}}5 $.   
\end{proposition}
\proof
Let $d$ be the dangerous center elected  by \D\ according to the rules of Phase 1.0. Let $a_1$ be the neighbor of $d$ defined in Lemma~\ref{lem:a1a2}.
Since \D\ plays according to the breaking-$P_5$ strategy, he plays $a_1$. This vertex and at least one of its white leaf neighbors turn from white to red. Moreover, the elected vertex $d$ turns from white to blue. 
Hence \D\ sizes at least 7 points, which in turn implies that the breaking-$P_5$ strategy is a greedy strategy. Applying Proposition~\ref{prop:buj}, we have that \go\ will end in at most $\displaystyle \displaystyle \frac{3n-e^{*}+c^{*}}5 $ turns. 
Moreover, any critical center created at turn $k>k^{**}$ must belong to $D_{k-1}$. 
But $D_k=\emptyset$,  for all $k\geq k^{**}$. Hence, by Lemma~\ref{lem:0}, we have $c^*=c^{**}$. Finally, it is obvious that $e^*\geq e^{**}$. In conclusion, if \D\ follows the breaking-$P_5$ strategy, \go\ ends in at most $\displaystyle\frac{3n-e^{**}+c^{**}}5 $ turns. 
\qed

\subsection{The proof for \go}\label{sec:pfmain}

We now prove that if $F$ is a \sff\ forest on $n$ vertices, then $\displaystyle \ggg(F)\leq\frac{3n}{5} $. For this purpose, we introduce two new processes for each turn of the game. These processes come in addition to the coloring and values defined in the greedy strategy and they only concern blue vertices.
When turning blue, some vertices will also be highlighted. We denote by $H_k$ the set of highlighted vertices at the end of turn $k\geq 1$.
The highlighting process will only occur at \D's turns. Moreover all blue vertices will get either weight $0$ or $1$. We denote by $w_k:B_k\rightarrow\{0,1\}$ the weight function at the end of turn $k\geq 1$. 
For $v\in V(F)$ and $X,Y\in \{B,R,W\}$, we write $v:X\rightsquigarrow Y$ to indicate that the color of $v$ changes from $X$ to $Y$ during the current turn. 

\medskip\noindent
\textbf{Highlighting Process.} 
\begin{itemize}
 \item When \D\ elects the dangerous center $d\in D_{k-1}$ at turn $k\geq 1$, $d$ is highlighted.    
\end{itemize}
Note that $H_k\subseteq B_k$, because the elected vertex always turns blue during the turn where it is elected. 

\medskip\noindent
\textbf{Weighting Process.}
\begin{itemize}
 \item \D's turn $k$, $k\geq 1$. All new blue vertices which are not highlighted get weight $1$. The new highlighted vertex $d$ gets weight $0$.
 \item \St's turn $k$, $k\geq 1$. She selects vertex $v$.
\begin{enumerate}[(A)]
 \item $v:W\rightsquigarrow R$ and $v$ is not the $P_5$-child of a vertex $d\in H_{k-1}$. All new blue vertices get weight $1$.
 \item $v:W\rightsquigarrow R$ and $v$ is the $P_5$-child of a vertex $d\in H_{k-1}$. All new blue vertices get weight $0$ and the vertex $d$ gets weight $1$.
 \item $v:B\rightsquigarrow R$, $v\not\in H_{k-1}$ and $v$ is not the $P_5$-child of a vertex $d\in H_{k-1}$. All new blue vertices get a weight equal to the weight of $v$.
\item $v:B\rightsquigarrow R$, $v\not\in H_{k-1}$ and $v$ is the $P_5$-child of a vertex $d\in H_{k-1}$. All new blue vertices get weight $0$ and $d$ get weight $1$.
 \item $v:B\rightsquigarrow R$ and $v\in H_{k-1}$. In particular $v$ is in $D_0$. Let $a_2$ be the $P_5$-child of $v$. 
                \begin{enumerate}
                       \item If $a_2\in W_{k-1}$, all new blue vertices but $a_2$ get weight $1$ and $a_2$ gets weight $0$.
                       \item If $a_2\in B_{k-1}$, all new blue vertices get weight $1$ and $a_2$ gets new weight $0$.
                       \item If $a_2\in R_{k-1}$, all new blue vertices get weight $1$.
                \end{enumerate}
\end{enumerate}
\item Additional process for turn $k$, $k\geq 1$. This process is applied after the two former processes. If the $P_5$-child of a vertex $d\in H_{k-1}$ turns from blue to red without being played by one of the player, then $d$ gets a weight of $1$.
\end{itemize}

Before being able to complete the proof, we need some technical definitions and lemmas. Until the end of this subsection, we always assume that \D\ plays according to the breaking-$P_5$ strategy.
\begin{definition}
Let $u$ be a vertex of a rooted forest $F$. For each turn $k\geq 0$, we define the following connected subtrees of $F\setminus R_k$. 
\begin{itemize}
\item $C_k(u)$ is the largest connected subtree of $F\setminus R_k$ which contains $u$.
\item $C^+_k(u)$ is the largest connected subtree of $F\setminus R_k$ which contains $u$ and whose vertices have height greater or equal to $h(u)$.
\item $C^-_k(u)=C_k(u)\setminus C^+_k(u)$.
\end{itemize}
\end{definition}
Remark that all these three sequences are non increasing with respect to inclusion. The following result follows directly from the definition.
\begin{lemma}\label{lem:1}
 Let $F$ be a rooted forest. For all $k\geq 0$ and all $u,v\in V(F)$, if $v$ is in $C^+_k(u)$ (resp. in $C^-_k(u)$), then $C^+_k(v)\subseteq C^+_k(u)$ 
(resp. $C^-_k(v)\subseteq C^-_k(u)$).  
\end{lemma}
\begin{lemma}\label{lem:2}
 Let $F$ be a rooted forest and $d$ be a vertex in $D_0(F)$, which is highlighted at turn $k\geq1$. Then $C^+_{k'}(d)\cap D_{k'}=\emptyset$, for all $k'\geq k$.
\end{lemma}
\proof
According to the prescribed strategy for \D, a dangerous center is highlighted at turn $k\geq1$ only if it has a maximum height among all the dangerous center from $D_{k-1}$. 
In other words, $C^+_{k-1}(d)\cap D_{k-1}=\{d\}$. Moreover, $C^+_{k'}(d)\subseteq C^+_{k-1}(d)$ and $D_{k'}\subseteq D_{k-1}$, for all $k'\geq k$. 
Hence,  $C^+_{k'}(d)\cap D_{k'}\subseteq \{d\}$. Finally, since $d$ belongs to $B_k$, it is not in $D_{k'}$. We conclude that  $C^+_{k'}(d)\cap D_{k'}=\emptyset$.
\qed

\begin{lemma}\label{lem:hlight} Let $F$ be a rooted forest. Let $d$ and $d'$ be two vertices in $D_0(F)$ which are respectively highlighted at turns $k$ and $k'\geq 1$. Then the following hold. 
 \begin{enumerate}[(i)]
    \item The $P_5$-child $a_2$ of the vertex $d$ belongs to $W_k$. Moreover, if $a_2$ is going to turn red or blue during Phase 1.0, this could be only during a turn of \St.
    \item The vertex $d$ is not the $P_5$-child of $d'$. 
  \end{enumerate}  
\end{lemma}
\proof 
We prove (i) first. Let $d$ be a dangerous center highlighted at turn $k$. The vertex $d$ has to belong to $D_{k-1}$. By definition of the sets of dangerous centers, this implies that its $P_5$-child $a_2$ is in $W_{k-1}$. 
Since $d$ is highlighted at turn $k$ this turn belongs to \D\ and during this turn he plays the vertex $a_1$ related to $d$. This vertex is not adjacent to the white vertex $a_2$. 
Hence $a_2$ remains white at the end of turn $k$. Moreover, by Lemma~\ref{lem:2}, we have $C^+_{k}(d)\cap D_k=\emptyset$. This implies that \D\ will not play 
in $C^+_{k}(d)$ until the end of Phase 1.0. In other words, \D\ will play all his remaining moves in Phase 1.0 in $C^-_{k}(d)$. 
But $a_2\in W_k$ is not adjacent to any vertex in $D^-_{k}(d)$.
Hence, \D's further moves in Phase 1.0 will not change the color of $a_2$. 

\medskip\noindent
We now prove (ii). Assume first that $k'>k$. By (i), the $P_5$-child of $d'$ belongs to $W_{k'}$. But $d$ has turned blue at turn $k$. 
Hence $d$ cannot be the $P_5$-child of $d'$. 
Suppose now that $k'<k$. By Lemma~\ref{lem:2}, $C^+_{k-1}(d')\cap D_{k-1}$ is empty. Therefore the $P_5$-child of $d'$ does not belong to $D_{k-1}$. But to be highlighted at turn $k$, $d$ has to belongs to $D_{k-1}$. We conclude as previously.  
\qed

We are now ready to state the following proposition, which will directly imply our main result for \go. 
\begin{proposition}\label{prop:main} Let $F$ be a \sff\ forest. For all $k\in\{0,...,k^{**}\}$, define $S_k=\displaystyle\sum_{i=1}^k e_i$ and $\displaystyle K_k=\sum_{v\in B_k} w_k(v)$. The following statements are true.
\begin{itemize}
 \item $I_k:$ $K_k\leq S_k$.
 \item $II_k:$ $\forall v\in B_k\setminus H_k$, if $w_k(v)=0$, then $C^-_k(v)=\emptyset$ and $C^+_k(v)\cap D_k=\emptyset$.
 \item $III_k:$  $\forall v\in B_k$, if $w_k(v)=0$, then $v$ is not a critical center at the end of turn $k$.
\end{itemize}
\end{proposition}
\proof
We will proceed by induction on $k$. Since $B_0$ is empty and $K_0=S_0=0$, the three statements are trivially true for $k=0$. Suppose that $I_k$, $II_k$ and $III_k$ are true for some $k\in\{0,...,k^{**}-1\}$. 
We denote by $b_{k+1}$ the number of vertices which turn blue during turn $k+1$.   

First, we deal with the weighting process for \D's turn. We assume that turn $k+1$ belongs to \D. At least two vertices turn from white to red and \D\ earns $6$ points with these two vertices. Hence, $e_{k+1}\geq b_{k+1}-1$, which implies that $S_{k+1}\geq S_k+b_{k+1}-1$. 
During this turn the weight of the vertices in $B_k$ does not change and only $b_{k+1}-1$ new blue vertices get weight $1$ (the highlighted vertex gets weight $0$). Therefore, $K_{k+1}\leq K_k+b_{k+1}-1$. 
Using the induction hypothesis, it yields $K_{k+1}\leq S_{k+1}$. We conclude that $I_{k+1}$ is true. 

\medskip\noindent
Let $w$ be a vertex in $B_k$, such that $w_k(w)=0$. Since $C^+_{k+1}(w)\cap D_{k+1}\subseteq C^+_{k}(w)\cap D_{k}$ and $C^-_{k+1}(w)\subseteq C^-_{k}(w)$, the truth 
of statement $II_{k+1}$ is straightforward from the induction hypothesis $II_k$ for the old blue vertices whose weight remains $0$. In other words, we only have to pay attention to blue vertices, new or old, which get new weight $0$ during turn $k+1$. 
Here, the weight of the former blue vertices does not change during this turn and all the new blue vertices which are not in $H_{k+1}$ get weight $1$. Hence $II_{k+1}$ is trivially true.
By Observation~\ref{ob:1}, in order to prove $III_{k+1}$, we only need to pay attention to vertices whose weight turns to $0$ during this turn. 
Here, it is the case for an unique vertex, the dangerous center  $d\in D_k$, which is elected by \D\ and highlighted during turn $k+1$. We emphasize that this is the only point of the whole proof where the \sff\ condition is used. 
Since \D\ follows the breaking-$P_5$ strategy, it implies that the vertex $a_1$ (defined in Lemma~\ref{lem:a1a2}) is played by \D\ at turn $k+1$. Hence, $a_1$ belongs to $R_{k+1}$. 
Let now $P$ be a path in $\mathcal P(F)$ whose center is $d$. By Lemma~\ref{lem:a1a2}, $P$ goes through $a_1$ which is in $R_k$. Therefore $P$ cannot be a critical $P_5$ in the residual forest $F_{k+1}$, which in turn implies 
that $d$ is not a critical center at turn $k+1$. In conclusion $III_{k+1}$ holds.

\medskip
Now, we prove the heredity of our statements for the weighting process of \St's turn. We suppose that turn $k+1$ belongs to \St\ and denote by $v$ the vertex played by her. We have to deal with the five different cases of the weighting process. 
We start by proving $I_{k+1}$. For the two first cases, \St\ earns $3$ points with the vertex $v$ which turns from white to red. Hence $e_{k+1}\geq b_{k+1}$ and $S_{k+1}\geq S_k+b_{k+1}$. 

\medskip\noindent
\textbf{Case (A).~}Because all new blue vertices get weight $1$ and the weight of former blue vertices does not change, we have $K_{k+1}\leq K_k+b_{k+1}$. 
By induction hypothesis, $K_k\leq S_k$, so we conclude that $W_{k+1}\leq S_{k+1}$.

\medskip\noindent
\textbf{Case (B).~}Since $v$ is the $P_5$-child of the vertex $d\in H_{k}$, it is adjacent to a leaf $l$ of $F$. Since $v$ belongs to $W_k$, this leaf is also in $W_k$. 
Hence, at least two vertices, $v$ and $l$ turn from white to red during turn $k+1$. It implies that $e_{k+1}\geq 3$ and $S_{k+1}\geq S_k+3$. 
Since all new blue vertices get weight $0$ and $d$ is the only former blue vertex whose weight is increased to $1$, we get $K_{k+1}\leq K_k+1\leq S_k+1 < S_{k+1}.$ 

\medskip
\noindent
For the three last cases, the vertex $v$ played by \St\ turns from blue to red and \St\ earns $2$ points with this vertex. Hence, $e_{k+1}\geq b_{k+1}-1$ and $S_{k+1}\geq S_k+b_{k+1}-1$.

\medskip\noindent
\textbf{Case (C).~}If $w_k(v)=0$, then all new blue vertices get weight $0$ and $K_{k+1}\leq K_{k}$. By induction, we get $K_{k+1}\leq S_k\leq S_{k+1}$.
Now, if $w_k(v)=1$, all new blue vertices get weight $1$. But $v$ will not be blue anymore, that is to say $v\not\in B_{k+1}$. Hence, its weight $1$ will not count for $K_{k+1}$.
Therefore, $K_{k+1}\leq K_k+b_{k+1}-1\leq S_{k+1}$.

\medskip\noindent
\textbf{Case (D).~}As above, the weight of $v$ will not count for $K_{k+1}$. Hence $K_{k+1}\leq K_k+1-w_k(v)$. Since $v$ is the $P_5$-child of a vertex $d\in H_k\subseteq B_k$, $C^-_k(v)$ contains at least $d$. 
By Lemma~\ref{lem:hlight} (ii), $v$ is not in $H_k$. Hence, the induction hypothesis $II_k$ yields $w_k(v)=1$. 
In conclusion, $K_{k+1}\leq K_k+1-w_k(v)\leq K_k\leq S_k\leq S_{k+1}$.

\medskip\noindent
\textbf{Case (E).~}Let $a_2$ be the $P_5$-child of $v$. 
\begin{enumerate}[(a)]
 \item If $a_2\in W_k$, then $K_{k+1}\leq K_k+b_{k+1}-1$, because $a_2$ gets weight $0$. By induction, $K_{k+1}\leq S_k+b_{k+1}-1\leq S_{k+1}$.
 \item If $a_2\in B_k$, $a_2$ gets new weight $0$. Therefore, $K_{k+1}\leq K_k+b_{k+1}-w_k(a_2)$. Since $v\in H_k\subseteq B_k$, $C^-_k(a_2)$ contains $v$. 
By Lemma~\ref{lem:hlight}, $a_2$ is not in $H_k$. Finally, $II_k$ implies that $w_k(a_2)=1$, which in turn implies $K_{k+1}\leq K_k+b_{k+1}-1\leq S_{k+1}$.
\item If $a_2\in R_k$, we have $K_{k+1}\leq K_k+b_{k+1}-w_k(v)$, because $v$ does not belong to $B_{k+1}$. By Lemma~\ref{lem:hlight}, we know that $a_2$ has turned red some moves after the moment $v$ has been highlighted. Moreover,
it happened during a \St's turn. There are two possibilities. First, $a_2$ has turned red because it has been played by \St\ at turn $k'<k+1$. 
Case~B or Case~D of the weighting process for \St's turn has been applied. 
Hence, during this turn, $v$ has gotten weight $1$. Second, $a_2$ might have turned red without being played by \St.
The additional weighting process ensures that $d$ has also gotten weight $1$. Moreover, the only weighting process, which can change the former weight of a blue vertex 
into $0$, is the second item of Case~E for \St's turn. But, this process is applied only if $v$ is the $P_5$-child of another highlighted vertex. 
Once more, Lemma~\ref{lem:hlight} ensures that it is impossible because $v$ is itself highlighted. In all the cases $w_k(v)=w_{k'}(v)=1$ and we can conclude as before.      
\end{enumerate}

\medskip\noindent
We now prove $II_{k+1}$ and $III_{k+1}$. As for \D's turn, we only have to pay attention to the blue vertices, whose weight becomes $0$ during turn $k+1$.

\medskip\noindent
\textbf{Case (A).} None of the new blue vertices gets weight $0$ and no former blue vertex weight change from $1$ to $0$. Then, $II_{k+1}$ and $III_{k+1}$ 
are trivially true.

\medskip\noindent
\textbf{Case (B).} The only vertices whose weight turn to $0$ during the turn are the new blue vertices in $N(v)$. Let $u$ be such a vertex. The vertex $v$ is the $P_5$-child of a vertex $d\in H_k$. 
Note that because $d$ was already blue at the beginning of turn $k+1$, we have $d\neq u$. 
The vertex $d$ is the only neighbor of $v$ whose height is strictly less than the height of $v$. Hence $h(u)>h(v)>h(d)$. Since $v$ is now red, $C^-_{k+1}(u)$ is empty. 
The vertex $d$ has been highlighted during a \D's turn, say turn $k'<k+1$. Since $v$ and $u$ are in $W_{k'}$, $u$ belongs to $C_{k'}^+(d)$.
By Lemma~\ref{lem:1}, $C^+_{k+1}(u)\cap D_{k+1}\subseteq C^+_{k'}(u)\cap D_{k'}\subseteq C^+_{k'}(d)\cap D_{k'}$.
By Lemma~\ref{lem:2}, the last set of these inclusion list is empty. Thus $C^+_{k+1}(u)\cap D_{k+1}$ is also empty. In conclusion $II_{k+1}$ is proved.
Moreover, we have $C^+_{k'}(d)\cap D_{k'}=\emptyset$ and $u\in C^+_{k'}(d)$. Hence $u$ is not in $D_{k'}$, but it is in $W_{k'}$. 
By Lemma~\ref{lem:0}, $u$ is not a critical center at the end of turn $k+1$. In conclusion, $III_{k+1}$ holds.

\medskip\noindent
\textbf{Case (C).} If $w_k(v)=1$, then no blue vertices, old or new, would see their weight turn to $0$. In that case $II_{k+1}$ and $III_{k+1}$ are trivially true.
Suppose that $w_k(v)=0$ and let $u\in N(v)$ be a new blue vertex. The vertex $v$ is not in $H_k$. By induction hypothesis $II_k$, $C^-_k(v)=\emptyset$ and $C^+_k(v)\cap D_k=\emptyset$. Hence, $u$ has to belong to $C^+_k(v)\setminus D_k$. 
By Lemma~\ref{lem:0}, $u$ cannot be a critical center at turn $k+1$. That proves $III_{k+1}$. By Lemma~\ref{lem:1}, $C^+_{k+1}(u)\subseteq C^+_{k+1}(v)$. Therefore $C^+_{k+1}(u)\cap D_{k+1}=\emptyset$. 
Finally $C^-_{k+1}(u)=\emptyset$, because $v$ is red at the end of turn $k+1$. In conclusion $II_{k+1}$ is true. 

\medskip\noindent
\textbf{Case (D).} The proof is the same as for Case (B).

\medskip\noindent
\textbf{Case (E).} Let $a_2$ be the $P_5$-child of $v$. For subcases (a) and (b), this vertex is the only blue vertex which gets weight $0$ during this turn. The vertex $a_2$ belongs to $C^+_k(v)$. 
It implies that $C^+_{k+1}(a_2)\subseteq C^+_{k}(a_2)\subseteq C^+_k(v)$ and $C^+_{k+1}(a_2)\cap D_{k+1}\subseteq C^+_{k}(v)\cap D_{k}$. Since $v$ is an highlighted vertex, $C^+_k(v)\cap D_k$ is empty by Lemma~\ref{lem:2}. 
Hence $C^+_{k+1}(a_2)\cap D_{k+1}$ is also empty. Moreover, because $v$ is red at the end of turn $k+1$, we have $C^-_{k+1}(a_2)=\emptyset$. That shows that $II_{k+1}$ holds.    
Since $v\in H_k$, this vertex has been highlighted during a \D's turn, say turn $k'\leq k$. We have $C^+_{k'}(v)\cap D_{k'}=\emptyset$. Hence, $a_2$ does not belongs to $D_{k'}$. 
By Lemma~\ref{lem:hlight}, $a_2$ belongs to $W_{k'}$. Therefore, applying Lemma~\ref{lem:0}, the vertex $a_2$ is not a critical center at the end of turn $k+1$.
In conclusion $III_{k+1}$ is true. Finally, for subcase (c), no blue vertex get new weight $0$, so there is nothing to prove.

\medskip
Finally, we have to prove that applying the additional weighting process at the end of turn $k+1$ does not change the true of statement $I_{k+1}$, $II_{k+1}$ and $III_{k+1}$. 
Since no blue vertex gets weight $0$ during this process, $II_{k+1}$ and $III_{k+1}$ remain trivially true. 
Let $x_1,...,x_m$ be the $P_5$-children of some vertices in $H_k$, which turn red without being played by one of the player. 
After applying this process, the new value of $K_{k+1}$ will be less or equal to $K_{k+1}+m-w_k(x_1)-\cdots -w_k(x_m)$. By the $P_5$-child definition, $C^-_k(x_i)\cap H_k$ is non empty, for any $i\in\{1,...,m\}$. Moreover, by Lemma~\ref{lem:hlight}, we have $x_i\not\in H_k$. Applying the induction hypothesis $II_k$, we get that $w_k(x_i)=1$.  
That proves that the value of $K_{k+1}$ is not increased by the application of the additional process. We conclude that $I_{k+1}$ remains true. \qed
\begin{theorem}
 If $F$ is a \sff\ forest on $n$ vertices, then $\displaystyle  \ggg(F)\leq\frac {3n}5$.
\end{theorem}
\proof
By Proposition~\ref{prop:1}, we only have to prove that $c^{**}\leq e^{**}$. 
Let $K_{k^{**}}$ and $S_{k^{**}}$ be defined as in Proposition~\ref{prop:main}. Note that $S_{k^{**}}=e^{**}$. Statement $III_{k^{**}}$ from Proposition~\ref{prop:main} implies that $c^{**}\leq K_{k^{**}}$. 
Moreover, statement $I_{k^{**}}$ implies that $K_{k^{**}}\leq S_{k^{**}}$. We conclude that $c^{**}\leq e^{**}$.  
\qed

\subsection{The proof for \gt}

The proof for \gt\ will proceed in the same way as for \go, except that now we introduce Phase 0 of the game, which corresponds to the first move of \St. Hence, Phase 1.1 will now start at turn $k=2$. 
We still write $k^{**}$ for the last turn of Phase 1.1 and $c^{**}$ for the number of critical centers at the end of turn $k^{**}$. 
We set $e_1=p(F_{0})-p(F_1)-5$. For $k\geq 2$, $e_i$ is defined as for \go. We also set $e^{**}=e_2+\cdots+e_{k^{**}}$.

In the proof of \cite[Theorem~2]{bujtas-2014}, it is proved that, by following the greedy strategy in \gt\ for an isolate-free forest on $n$ vertices, \D\ can force the game to end in at most $\displaystyle  \frac{3n-e_1-e^*+c^*} 5 $ turns. 
Here, $e^*$ is the amount of bonus point earned during Phase 1 and $c^*$ is the number of critical centers at the end of this same phase. As for \go, this statement implies the following proposition.

\begin{proposition}
 Let $F$ be an isolate-free forest on $n$ vertices. If \D\ plays \gt\ according to the breaking-$P_5$ strategy, then the number of turns will be at most $\displaystyle \frac{3n-e_1-e^{**}+c^{**}}5 $.   
\end{proposition}
We can define the weighting and the highlighting process exactly in the same manner as for \go. We can state a proposition very similar to Proposition~\ref{prop:main}.
Only statement $I_k$ will change. For this new proposition, we will have $K_k\leq S_k+2$, for all $k\geq 0$. Indeed $e_1$ is equal to the number of blue vertices minus $2$. Since all these blue vertices will get weight $1$, we have $K_1\leq S_1+2$.
Finally, this new proposition will directly prove that $c^{**}\leq e_1+e^{**}+2$. It clearly implies the below theorem.
\begin{theorem}
 If $F$ is a \sff\ forest on $n$ vertices, then $\displaystyle  \gggp(F)\leq\frac {3n+2}5$.
\end{theorem}

\section*{Acknowledgements}

The author would like to thank Sandi Klav\v zar and Csilla Bujt\'as for helpful discussions about this work. In particular, the nice definition of \sff\ forest owes a great deal to Sandi Klav\v zar. 
The research was financed by the ANR-14-CE25-0006 project of the French National Research Agency and by the grant CMIRA Explo'RA Doc from La R\'{e}gion Rh\^{o}ne Alpes.


\begin{thebibliography}{}
\bibitem{complexity-2014+}
 B.~Bre{\v{s}}ar, P.~Dorbec, S.~Klav{\v{z}}ar, G.~Ko\v smrlj,
 Complexity of the game domination problem, 
 manuscript, 2014. 

\bibitem{dorbec-2015+}
 B.~Bre{\v{s}}ar, P.~Dorbec, S.~Klav{\v{z}}ar, G.~Ko\v smrlj,
 How long one can bluff in the domination game?, 
 manuscript, 2015. 


\bibitem{bresar-2010}
  B.~Bre{\v{s}}ar, S.~Klav{\v{z}}ar, D.~F.~Rall,
  Domination game and an imagination strategy,
  SIAM J.\ Discrete Math.\ 24 (2010) 979--991.

\bibitem{bresar-2013}
  B.~Bre{\v{s}}ar, S.~Klav{\v{z}}ar, G.~Ko\v smrlj, D.~F.~Rall,
  Domination game: extremal families of graphs for the 3/5-conjectures,
  Discrete Appl.\ Math.\ 161 (2013) 1308--1316.

\bibitem{brdo-2014}
 B.~Bre{\v{s}}ar, P.~Dorbec, S.~Klav{\v{z}}ar, G.~Ko\v smrlj,
 Domination game: Effect of edge- and vertex-removal,
 Discrete Math.\ 330 (2014) 1--10.

\bibitem{brkl-2013b}
  B.~Bre{\v{s}}ar, S.~Klav{\v{z}}ar, D.~F.~Ral,
  Domination game played on trees and spanning subgraphs,
  Discrete Math. 313 (2013) 915--923.

\bibitem{bujtas-2013}
  Cs.~Bujt\'as,
  Domination game on trees without leaves at distance four,
  Proceedings of the 8th Japanese-Hungarian Symposium on Discrete Mathematics and Its Applications (A.~Frank, A.~Recski, G.~Wiener, eds.) 73--78, (June 2013).

\bibitem{bujtas-2015}
  Cs.~Bujt\'as,
  On the game domination number of graphs with given minimum degree,
  Electron.\ J.\ Combin., to appear.

\bibitem{bujtas-2014}
  Cs.~Bujt\'as, Domination game on forests, Discrete Math. 338 (2015) 2220--2228. 


\bibitem{dorbec-2015}
  P.~Dorbec, G.~Ko\v smrlj, G.~Renault,
  The domination game played on unions of graphs,
  Discrete Math.\ 338 (2015) 71--79.
 
\bibitem{henning-2014+}
  M.~A.~Henning, W.~B.~Kinnersley,
  Bounds on the game domination number,
  manuscript, 2014.

\bibitem{henning-2015}
  M.~A.~Henning, S.~Klav{\v{z}}ar, D.~F.~Rall,
  Total version of the domination game,
  Graphs Combin., in press, DOI 10.1007/s00373-014-1470-9.

\bibitem{henning-2015+}
  M.~A.~Henning, S. Klav\v zar, D.~F.~Rall,
  The 4/5 upper bound on the game total domination number,
  Combinatorica, to appear. 

\bibitem{kinnersley-2013}
  W.~B.~Kinnersley, D.~B.~West, R.~Zamani,
  Extremal problems for game domination number,
  SIAM J.\ Discrete Math.\ 27 (2013) 2090--2107.

\bibitem{kosmrlj-2014}
  G.~Ko\v{s}mrlj,
  Realizations of the game domination number,
  J. Comb.\ Optim.\ 28 (2014) 447--461.

\bibitem{javad-2015+}
  M.~J.~Nadjafi-Arani, M.~Siggers, H.~Soltani,
  Characterisation of forests with trivial game domination numbers,
  J. Comb.\ Optim., in press, DOI 10.1007/s10878-015-9903-9.  

\end{thebibliography}
\end{document}